\documentclass[a4paper,12pt]{article}
\usepackage{CJK,CJKnumb}
\usepackage{indentfirst}
\usepackage{amssymb}
\usepackage{amsmath}
\usepackage[dvips]{graphicx}
\usepackage{amsthm}
\usepackage{amscd}
\usepackage{latexsym}
\usepackage{mathrsfs}

    \date{}
    \allowdisplaybreaks[4] \footskip=15pt

    \numberwithin{equation}{section} \theoremstyle{plain}
    \newtheorem*{thm*}{Main Theorem}

    \newtheorem*{cor*}{Corollary}
    
    \newtheorem*{lem*}{Lemma}
    
    \newtheorem*{prop*}{Proposition}
    
    \newtheorem*{rem*}{Remark}
    
    \newtheorem*{exa*}{Example}
    
    \newtheorem*{df*}{Definition}
    
    \newtheorem*{conj*}{Conjecture}
    \newtheorem*{ack*}{ACKNOWLEDGEMENTS}

    \newcommand{\pf}{\noindent\begin {proof}}
    \newcommand{\epf}{\end{proof}}

\begin{document}
\begin{center}
{\Large\bf Prime ideals in decomposable lattices}\vskip 8mm

{\bf Xinmin Lu$^{*}$\,\,\,\,\,\,\, Dongsheng Liu} \vskip 1mm

{\small\it School of Science, Nanjing University of Science and
Technology,\\ Nanjing 210094, P.R. China}\vskip 2mm

{\bf Zhinan Qi} \vskip 1mm

{\small\it Department of Mathematics, Nanchang University,\\
Nanchang 330047, P.R. China}\vskip 2mm

{\bf Hourong Qin} \vskip 1mm

{\small\it Department of Mathematics, Nanjing University,\\
Nanjing 210093, P.R. China}\vskip 4mm

\begin{minipage}{130mm}
{\noindent\small\bf Abstract.} A distributive lattice $L$ with
minimum element $0$ is called decomposable lattice if $a$ and $b$ are not
comparable elements in $L$ there
exist $\overline{a},\overline{b}\in L$ such that
$a=\overline{a}\vee(a\wedge b), b=\overline{b}\vee(a\wedge b)$ and
$\overline{a}\wedge \overline{b}=0$. The main purpose of this
paper is to investigate prime ideals, minimal prime ideals and
special ideals of a decomposable lattice. These are keys to understand the algebraic structure of decomposable
lattices. \vskip 2mm

\vskip 2mm

\vskip 2mm
\end{minipage}
\end{center}
$\\$

\footnotetext {$^{*}$ Corresponding author: School of Science,
Nanjing University of Science and Technology, Nanjing 210094, P.R.
China. {\bf E-mail}: xmlu\_nanjing@hotmail.com}

\noindent{\bf 1. Introduction }\vskip 4mm

In [10] Gr\"{a}tzer and Schmidt characterized a Stone lattice as a
distributive pseudocomplemented lattice in which every prime ideal
contains a unique minimal prime ideal. Motivated by this
characterization of Stone lattices, Cornish and Pawar
characterized distributive lattices with minimum element $0$ in
which each prime ideal contains a unique minimal prime ideal (see
e.g. [4,13]) and distributive lattices with $0$ in which each
prime ideal contains $n$ minimal prime ideals [5]. They called
such lattices respectively normal lattices and $n$-normal
lattices. As a natural generalization of normal lattices, Cornish
also introduced the concept of relatively normal lattices: a
relatively normal lattice is a distributive lattice with $0$ such
that every bound closed interval is a normal lattice [5]. Filipoiu
and Georgescu investigated values (regular ideals) in relatively
normal lattices [7]. Hart, Snodgrass and Tsinakis further studied
the structure of relatively normal lattices (see e.g. [11,14,15]).
Motivated by the above works, we shall be concerned with
decomposable lattices by replacing the "normality" by the
"decomposability", i.e., a decomposable lattice is a distributive
lattice $L$ with minimum element $0$ such that for any $a,b\in L$,
if  $a$ and $b$ are not comparable elements, written by
$a\parallel b$, then there exist $\overline{a},\overline{b}\in L$
such that $a=\overline{a}\vee(a\wedge b),
b=\overline{b}\vee(a\wedge b)$ and $\overline{a}\wedge
\overline{b}=0$.

Decomposability is not just the algebraic properties for some
lattices. There exist in other algebraic areas, such as rings,
modules and lattice-ordered group. We will see examples in section
2. Decomposable lattice is the common tool to understand these
properties. Furthermore, the characterizations of prime ideals,
minimal prime ideals and special ideals in the decomposable
lattice are explicit. More details will be seen in later.
Moreover, these characterizations can be our main technical tool
for the further study of the structure of such lattices. In fact,
with the help of the results of the present paper, the structure
of decomposable lattices determined by their prime ideals, minimal
prime ideals and special ideals can be developed [12].

Here is a brief outline of the article. In Section 1, we simply
review some basic definitions and some well-known results.  Three
examples of decomposable lattices in lattices, rings and
lattice-ordered groups, respectively are given. In Section 2, we
investigate prime ideals of a decomposable lattice and the
relationship between prime ideals and regular ideals. This is
contained in Section 3, where we shall first establish explicit
characterizations of minimal prime ideals of a decomposable
lattice and then investigate the relationship among prime ideals,
minimal prime ideals and regular ideals.  We investigate  special
ideals of a decomposable lattice and the relationship between
special ideals and regular ideals in the last section.\vskip 8mm

\noindent{\bf 2. Preliminaries and Examples}\vskip 4mm

Firstly, we simply review some basic definitions and some
well-known results. The reader is refereed to [9] for the general
theory of lattices.

Throughout this paper, we consider lattices $L$ with minimum
element $0$, denote by $\mathbb{DL}$ the class of decomposable
lattices and use "$\subset$" and "$\supset$" to denote proper
set-inclusion.

A lattice $L$ is called distributive if $a\wedge(b\vee c)=(a\wedge
b)\vee(a\wedge c)$ for any $a,b,c\in L$. A nonempty subset $I$ in
a lattice $L$ is called an ideal of $L$ if $a\vee b\in I$ for any
$a,b\in I$ and $a\geq x\in L$ implies that $x\in I$. We denote by
$Ide(L)$ the set of all ideals of $L$. In particular, if $a\in L$
then $(a]=\{x\in L|\,\, x\leq a\}$ is called the principal ideal
of $L$ generated by $a$. A direct computation shows that if $L\in
\mathbb{DL}$ then $Ide(L)$ is a distributive lattice by the rule:
$I\wedge J=I\cap J$ and $I\vee J=\{a\vee b|\,\, a\in I, b\in J\}$
for any $I,J\in Ide(L)$.

An ideal $P$ in a lattice $L$ is called prime if $P\neq L$ and
$a\wedge b\in P$ implies that either $a\in P$ or $b\in P$, where
$a,b\in L$. By Zorn's Lemma, each prime ideal contains a minimal
prime ideal. We denote by $Spe(L)$ and $MinSpe(L)$ respectively
the set of all prime ideals of $L$ and the set of all minimal
prime ideals of $L$.

Let $L$ be a lattice. For any $0<x\in L$, by Zorn's Lemma, there
exists a maximal ideal of $L$ with respect to not containing $x$,
denoted $M$, $M$ is called a regular ideal and is the value of
$x$. In general, $a$ need not have a unique value. We denote by
$Val(x)$ the set of all values of $x$. If $M$ is the unique value
of $x$, $M$ or $x$ is called special. We denote by $V(L)$ and
$S(L)$ respectively the set of all values of $L$ and the set of
all special values of $L$. Clearly, $S(L)\subseteq V(L)$. Observe
that the following conditions are equivalent: (1) $M\in V(L)$; (2)
$M$ is meet-irreducible, i.e., if
$\bigcap\limits_{\lambda\in\Lambda}I_{\lambda}=M$, where
$\{I_{\lambda}\}_{\lambda\in\Lambda}\subseteq Ide(L)$, then
$I_{\lambda}=M$ for some $\lambda$; (3) $M\subset
M^{*}=\bigcap\{I\in Ide(L)|\,\, I\supset M\}$; (4) $M\in Val(x)$,
where $x\in M^{*}\setminus M$.

For a lattice $L$ and $\emptyset\neq A\subseteq L$, we write
$A^{\perp}=\{x\in L|\,\, x\wedge a=0$ for any $a\in L \}$.
$A^{\perp}$ is called the polar of $A$, and define
$(A^{\perp})^{\perp}=A^{\perp\perp}$. $P\in Ide(L)$ is called
polar if $P=A^{\perp}$ for some $\emptyset\neq A\subseteq L$.
Clearly, $P\in Ide(L)$ is polar if and only if $P=P^{\perp\perp}$.
We denote by $P(L)$ the set of all polar ideals of $L$.

Let $L$ be a lattice. A nonempty subset $F$ of $L$ is called a
filter of $L$ if the following conditions are satisfied: (1)
$0\not\in F$; (2) for any $a,b\in F$, $a\wedge b\in F$; (3) if
$x\in L$ and $x\geq a\in F$ implies $x\in F$. By Zorn's Lemma,
each filter $F$ of $L$ must be contained in a maximal filter $U$
of $L$, and $U$ is called an ultrafilter of $L$.

We give the definition of decomposable lattice as following.\vskip
2mm

\noindent{\bf Definition 2.1. }A decomposable lattice is a
distributive lattice $L$ with minimum element $0$ such that for
any $a,b\in L$, if $a\parallel b$ then there exist
$\overline{a},\overline{b}\in L$ such that
$a=\overline{a}\vee(a\wedge b), b=\overline{b}\vee(a\wedge b)$ and
$\overline{a}\wedge \overline{b}=0$.\vskip 2mm

Followings are  examples of decomposable
lattices, which are closely related to rings and lattice-ordered
groups as well as lattices. \vskip 2mm

Recall that a lattice $L$ is called strongly projectable if $L=(a]
\vee a^{\perp}$ for any $a\in L$. \vskip 2mm

\noindent{\bf Example 2.2. }Let $L$ be a distributive lattice. If
$L$ is strongly projectable then $L\in \mathbb{DL}$.\vskip 2mm

\noindent{\bf Proof. }Given any $a,b\in L$ with $a\parallel b$,
since $L=(a\wedge b]\vee (a\wedge b)^{\perp}$, there exist
$x_{1},x_{2}\in (a\wedge b]$ and $y_{1},y_{2}\in (a\wedge
b)^{\perp}$ such that\vskip 2mm
\begin{center}
$a=x_{1}\vee y_{1},\,\,\, b=x_{2}\vee y_{2}$.\vskip 2mm
\end{center}
\noindent Since $L$ is distributive, we have\vskip 2mm
\begin{center}
$a\wedge b=(a\wedge b)\wedge (x_{1}\vee y_{1})=(a\wedge b\wedge
x_{1})\vee (a\wedge b\wedge y_{1})=a\wedge b\wedge x_{1}$.\vskip
2mm
\end{center}
\noindent So $a\wedge b\leq x_{1}$, which implies that
$x_{1}=a\wedge b$. Similarly, $x_{2}=a\wedge b$. Then\vskip 2mm
\begin{center}
$a=a\wedge(x_{1}\vee y_{1})=(a\wedge x_{1})\vee (a\wedge
y_{1})=(a\wedge y_{1})\wedge (a\wedge b)$.\vskip 2mm
\end{center}
\noindent Similarly, $b=(b\wedge y_{2})\wedge (a\wedge b)$.
Since\vskip 2mm
\begin{center}
$(a\wedge y_{1})\wedge (b\wedge y_{2})=(a\wedge b)\wedge
y_{1}\wedge y_{2}=0$, \vskip 2mm
\end{center}
\noindent we get $L\in \mathbb{DL}$.\vskip 2mm

Recall from [2,6] that a partially ordered group is both a group
$(G,+)$ and a partially ordered set $(G,\leq)$ whenever $a\leq b$
and $x,y\in G$ then $x+a+y\leq x+b+y$. A lattice-ordered group is
a partially ordered group $G$ and the underlying order is a
lattice. A lattice-ordered group is called complete if every
subset bounded above has a least upper bound and every subset
bounded below has a greatest lower bound [3]. Recall also from [1]
that a lattice-ordered group $G$ is called compactly generated if
$\{a_{\lambda}\}_{\lambda\in\Lambda}$ is a nonempty subset of $L$
and $\bigwedge\limits_{\lambda\in\Lambda}a_{\lambda}=0$ then there
exists a finite subset $\{a_{i}\}^{n}_{i=1}$ of
$\{a_{\lambda}\}_{\lambda\in\Lambda}$ such that
$\bigwedge\limits_{i=1}^{n}a_{i}=0$.\vskip 2mm

\noindent{\bf Example 2.3. }Let $(G,+,\vee,\wedge)$ be a complete
lattice-ordered group. If $G$ is compactly generated then the
positive cone $G^{+}=\{x\in G|\,\, x\geq 0\}\in
\mathbb{DL}$.\vskip 2mm

\noindent{\bf Proof. }By hypothesis, each positive element in $G$
can be written as a join of some atoms in $G$. So, for any $x,y\in
G^{+}$ with $x\parallel y$, write \vskip 2mm
\begin{center}
$x=\bigvee\limits_{\lambda\in \Lambda_{1}} a_{\lambda},\,\,\,
y=\bigvee\limits_{\mu\in \Lambda_{2}} b_{\mu}$, \vskip 2mm
\end{center}
\noindent where each $a_{\lambda}$ and $b_{\mu}$ are atoms in $G$.
If $\Lambda=\Lambda_{1}\cap \Lambda_{2}$, then $x\wedge
y=\bigvee\limits_{\nu\in \Lambda} c_{\nu}$. Now, set\vskip 2mm
\begin{center}
$x^{'}=\bigvee\limits_{\lambda\in \Lambda_{1}\setminus\Lambda}
a_{\lambda},\,\,\, y^{'}=\bigvee\limits_{\mu\in
\Lambda_{2}\setminus\Lambda} b_{\mu}$.\vskip 2mm
\end{center}
\noindent Then $x=x^{'}\vee(x\wedge y), y=y^{'}\vee(x\wedge y)$.
In view of [6], $G$ is completely distributive, we further have
\vskip 3mm
\begin{center}
$x^{'}\wedge y^{'}=(\bigvee\limits_{\lambda\in
\Lambda_{1}\setminus\Lambda}
a_{\lambda})\wedge(\bigvee\limits_{\mu\in
\Lambda_{2}\setminus\Lambda} b_{\mu})=\bigvee\limits_{\lambda\in
\Lambda_{1}\setminus\Lambda}\,\,\bigvee\limits_{\mu\in
\Lambda_{2}\setminus\Lambda}(a_{\lambda}\wedge b_{\mu})=0$.\vskip
3mm
\end{center}
\noindent So $G^{+}\in \mathbb{DL}$.\vskip 2mm

Following Fuchs [8], a ring $R$ is called arithmetical if the
lattice $Ide(R)$ of all ideals in $R$ is distributive, i.e.,
$I\cap(J+K)=(I\cap J)+(I\cap K)$ for any $I,J,K\in Ide(R)$.\vskip
2mm

\noindent{\bf Example 2.4. }If $R$ is an arithmetical ring and
satisfies that for any $I\in Ide(R)$ there exists some $e^{2}=e\in
R$ such that $I=eR$, then $Ide(R)\in \mathbb{DL}$.\vskip 2mm

\noindent{\bf Proof. }Given any $I,J\in Ide(R)$, if $I\parallel
J$, write $K=I\cap J\in Ide(R)$, then there exists some
$e^{2}=e\in R$ such that $K=eR$. Since $I\subseteq R=eR\oplus
(1-e)R$, there exist $I_{1},I_{2}\in Ide(R)$ with $I_{1}\subseteq
eR, I_{2}\subseteq (1-e)R$ such that $I=I_{1}+I_{2}$. Similarly,
there exist $J_{1},J_{2}\in Ide(R)$ with $J_{1}\subseteq eR,
J_{2}\subseteq (1-e)R$ such that $J=J_{1}+J_{2}$. Thus, we
have\vskip 2mm
\begin{center}
$K=K\cap I=K\cap (I_{1}+I_{2})=(K\cap I_{1})+(K\cap I_{2})=K\cap
I_{1}$,\vskip 2mm
\end{center}
\noindent and hence $K\subseteq I_{1}$, so that $K=I_{1}$.
Similarly, $K=J_{1}$. So\vskip 2mm
\begin{center}
$I=I\cap (I_{1}+I_{2})=(I\cap I_{1})+(I\cap I_{2})=(I\cap
I_{2})+(I\cap J)$\vskip 2mm
\end{center}
\noindent and \vskip 2mm
\begin{center}
$J=J\cap (J_{1}+J_{2})=(J\cap J_{1})+(J\cap J_{2})=(J\cap
J_{2})+(I\cap J)$.\vskip 2mm
\end{center}
\noindent Write $I^{'}=I\cap I_{2}, J^{'}=J\cap J_{2}$. Then
$I=I^{'}+(I\cap J), J=J^{'}+(I\cap J)$ and \vskip 2mm
\begin{center}
$I^{'}\cap J^{'}=(I\cap I_{2})\cap (J\cap J_{2})=(I\cap J)\cap
(I_{2}\cap J_{2})\subseteq eR\cap (1-e)R=0$.\vskip 2mm
\end{center}
\noindent Therefore $Ide(R)\in \mathbb{DL}$.\vskip 8mm

\noindent{\bf 3. Prime ideals}\vskip 4mm

In this section, we shall first establish characterizations of
prime ideals of a decomposable lattice and then investigate the
relationship between prime ideals and regular ideals.\vskip 2mm

\noindent{\bf Theorem 3.1. }Let $L\in \mathbb{DL}$ and $L\neq P\in
Ide(L)$. The following conditions are equivalent:

(1) $P\in Spe(L)$.

(2) If $x\wedge y=0$ then either $x\in P$ or $y\in P$ for $x,y\in
L$.

(3) $x,y\in L\setminus P$ implies $x\wedge y\in L\setminus P$.

(4) If $I\cap J\subseteq P$ then either $I\subseteq P$ or
$J\subseteq P$ for $I,J\in Ide(L)$.

(5) If $I,J\in Ide(L)$ such that $P\subseteq I$ and $P\subseteq
J$, then either $I\subseteq J$ or $J\subseteq I$. \vskip 2mm

\noindent{\bf Proof. }(1)$\Leftrightarrow$(2)$\Leftrightarrow$(3)
is clear.

(1)$\Rightarrow$(4) Let $I,J\in Ide(L)$ be such that $I\cap
J\subseteq P$. If $I\not\subseteq P$, then we may pick $a\in
I\setminus P$. So, for any $b\in J$, since $a\wedge b\in I\cap
J\subseteq P$ and $a\not\in P$, we get that $b\in P$, so that
$J\subseteq P$.

(4)$\Rightarrow$(1) Given $a,b\in L$ with $a\parallel b$, if
$a\wedge b\in P$, then $(a\wedge b]\subseteq P$. Notice that
$(a\wedge b]=(a]\cap(b]$. Then $(a]\cap(b]\subseteq P$. So, by
(4), we get that either $(a]\subseteq P$ or $(b]\subseteq P$,
hence either $a\in P$ or $b\in P$. Therefore $P\in Spe(L)$.

(1)$\Rightarrow$(5) Let $I,J\in Ide(L)$ be such that $P\subseteq
I$ and $P\subseteq J$. Suppose that $I$ and $J$ are not
comparable, written $I\parallel J$. Pick $a\in I\setminus J, b\in
J\setminus I$. Clearly, $a\parallel b$. Then there exist
$\overline{a},\overline{b}\in L$ such that
$a=\overline{a}\vee(a\wedge b), b=\overline{b}\vee(a\wedge b)$ and
$\overline{a}\wedge \overline{b}=0\in P$. So either
$\overline{a}\in P$ or $\overline{b}\in P$, which implies that
either $a\in J$ or $b\in I$, a contradiction.

(5)$\Rightarrow$(1) Given $a,b\in L$ with $a\wedge b\in P$, then
$(a]\vee P, (b]\vee P\supseteq P$. By (5), $(a]\vee P$ and
$(b]\vee P$ are comparable. Without loss of generality, assume
that $(a]\vee P\subseteq(b]\vee P$. Since $Ide(L)$ is a
distributive lattice, we then have\vskip 2mm
\begin{center}
$P=(a\wedge b]\vee P=(( a]\cap(b])\vee P=((a]\vee P)\cap((b]\vee
P)=(a]\vee P$.\vskip 2mm
\end{center}
\noindent Hence $a\in P$. Therefore $P\in
Spe(L)$.$\hfill\Box$\vskip 2mm

By Theorem 3.1, we now get some immediate corollaries which should
demonstrate some of the importance of prime ideals.\vskip 2mm

\noindent{\bf Corollary 3.2. }Let $L\in \mathbb{DL}$.

(1) $V(L)\subseteq Spe(L)$.

(2) $\bigcap V(L)=\bigcap Spe(L)=0$.

(3) For any $I\in Ide(L)$, $I=\bigcap\{M\in V(L)|\,\, M\supseteq
I\}=\bigcap\{P\in Spe(L)|\,\, P\supseteq I\}$.\vskip 2mm

\noindent{\bf Corollary 3.3. }Let $L\in \mathbb{DL}$.

(1) The intersection of a chain of prime ideals of $L$ is prime.

(2) If $P\in Spe(L)$ then the set $\{I\in Ide(L)|\, I\supseteq
P\}$ forms a chain.

(3) $L$ is totally ordered if and only if the zero ideal $0$ of
$L$ is prime.\vskip 2mm

\noindent{\bf Corollary 3.4. }Let $L\in \mathbb{DL}$. The
following conditions are equivalent:

(1) Each prime ideal of $L$ contains a unique minimal prime ideal.

(2) For any $N_{1},N_{2}\in MinSpe(L)$, if $N_{1}\parallel N_{2}$
then $L=N_{1}\vee N_{2}$.

(3) For any $P_{1},P_{2}\in Spe(L)$, if $P_{1}\parallel P_{2}$
then $L=P_{1}\vee P_{2}$.

(4) For any $M_{1},M_{2}\in V(L)$, if $M_{1}\parallel M_{2}$ then
$L=M_{1}\vee M_{2}$.\vskip 2mm

As an application of Theorem 3.1, we now investigate the
relationship between prime ideals and regular ideals of a
decomposable lattice.\vskip 2mm

\noindent{\bf Theorem 3.5. }Let $L\in \mathbb{DL}$. The following
conditions are equivalent:

(1) $Spe(L)=V(L)$.

(2) $V(L)$ satisfies $DCC$.

(3) $Spe(L)$ satisfies $DCC$ . \vskip 2mm

\noindent{\bf Proof. }(2)$\Leftrightarrow$(3) is clear. It
suffices to show (1)$\Leftrightarrow$(2)

(1)$\Rightarrow$(2) Given any descending chain of $V(L)$:
$Q_{1}\supseteq Q_{2}\supseteq\cdots\supseteq
Q_{n}\supseteq\cdots$, and set
$P=\bigcap\limits_{i=1}^{\infty}Q_{i}$. A direct computation shows
that $P\in Spe(L)$, hence $P\in V(L)$ by hypothesis. Then there
exists some positive integer $m$ such that $P=Q_{m}$. So
$Q_{1}\supseteq Q_{2}\supseteq\cdots\supseteq
Q_{m}=Q_{m+1}=\cdots$. Thus $V(L)$ satisfies $DCC$.

(2)$\Rightarrow$(1) Given any $P\in Spe(L)$, then\vskip 2mm
\begin{center}
$P=\bigcap\{Q\in V(L)|\,\, Q\supseteq P\}$.\vskip 2mm
\end{center}
\noindent Since $P$ is prime, by Corollary 3.3, the set $\{Q\in
V(L)|\,\, Q\supseteq P\}$ is a chain. Since $V(L)$ satisfies
$DCC$, this chain is finite, and hence it must have a least
element, denoted $Q_{0}$, so that $P=Q_{0}\in V(L)$. So
$Spe(L)=V(L)$. $\hfill\Box$\vskip 8mm

\noindent{\bf 4. Minimal prime ideals}\vskip 4mm

In this section, we first investigate the relationship between
ultrafilters and minimal prime ideals in a decomposable lattice.
With the help of the relationship, we shall establish explicit
characterizations of minimal prime ideals of a decomposable
lattice, which are pure lattice-theoretic extension of the
corresponding results of lattice-ordered groups [2,6].

Filters arise naturally whenever we have a partially ordered set.
We remind the reader that if $L$ is a lattice and $E$ is a
$\wedge$-semilattice of $L$ (i.e. for any $a,b\in E$, $a\wedge
b\in E$) then $\overline{E}=\bigcap\{F|\,\, E\subseteq F$ a filter
of $L\}$ is the smallest filter of $L$ containing $E$, is called
the filter of $L$ generated by $E$. By Zorn's Lemma, each filter
$F$ of $L$ must be contained in a maximal filter $U$ of $L$, and
$U$ is called an ultrafilter of $L$.\vskip 2mm

\noindent{\bf Lemma 4.1. }Let $L\in \mathbb{DL}$ and $U$ a
$\wedge$-semilattice with $0\not\in U$. The following conditions
are equivalent:

(1) $U$ is an ultrafilter of $L$.

(2) For any $x\in L\setminus U$ there exists $u\in U$ such that
$x\wedge u=0$.

(3) $L\setminus U\in MinSpe(L)$.\vskip 2mm

\noindent{\bf Proof. }(1)$\Rightarrow$(2) Assume that there exists
some $x\in L\setminus U$ such that for any $u\in U$, $x\wedge
u>0$. A direct computation shows that the set\vskip 2mm
\begin{center}
$U_{0}=\{x\wedge u|\,\, u\in U\cup\{x\}\}$\vskip 2mm
\end{center}
\noindent is a $\wedge$-semilattice with $0\not\in U_{0}$. Let
$\overline{U}$ be the filter of $L$ generated by $U_{0}$. Then
$U_{0}\subseteq\overline{U}$, so that $U\subseteq \overline{U}$.
But $x\in \overline{U}$ and $x\not\in U$, which contradicts the
fact that $U$ is an ultrafilter of $L$.

(2)$\Rightarrow$(3) We first show that $L\setminus U\in Ide(L)$.
Since $U$ is a filter, it suffices to show that for any $x,y\in
L\setminus U$, $x\vee y\in L\setminus U$. Assume that there exist
$x,y\in L\setminus U$ such that $x\vee y\in U$. By (2), there
exist $a,b\in U$ such that $x\wedge a=0, y\wedge b=0$. Set
$c=a\wedge b\in U$. Then $0=c\wedge(x\vee y)\in U$, a
contradiction.

We next show that $L\setminus U\in MinSpe(L)$. Clearly,
$L\setminus U\in Spe(L)$. Assume that $L\setminus U\not\in
MinSpe(L)$. Then by Zorn's Lemma, there exists some $M\in
MinSpe(L)$ such that $L\setminus U\supset M$. Observe that
$L\setminus M$ is a filter of $L$ and clearly $L\setminus M\supset
U$, which contradicts the fact that $U$ is an ultrafilter of $L$.
So $L\setminus U\in MinSpe(L)$.

(3)$\Rightarrow$(1) Clearly, if $L\setminus U\in MinSpe(L)$ then
$U$ is a filter of $L$. Assume that $U$ is not an ultrafilter of
$L$. Then there exists an ultrafilter of $L$, denoted $W$, such
that $W\supset U$. Using the result of (1)$\Rightarrow$(2),
$L\setminus W\in MinSpe(L)$. But $L\setminus W \supset L\setminus
U\in MinSpe(L)$, a contradiction. Therefore $U$ is an ultrafilter
of $L$.$\hfill\Box$\vskip 2mm

\noindent{\bf Lemma 4.2. }Let $L\in \mathbb{DL}$. If $X$ is a
$\wedge$-semilattice with $0\not\in X$ then\vskip 2mm
\begin{center}
$\bigcup\{a^{\perp}|\,\, a\in X\}=\bigcap\{P\in Spe(L)|\,\, P\cap
X=\emptyset\}=\bigcap\{M\in MinSpe(L)|\,\, M\cap X=\emptyset\}$.
\vskip 2mm
\end{center}
\noindent In particular, if $P\in Spe(L)$ then\vskip 2mm
\begin{center}
$\bigcup\{a^{\perp}|\,\, a\in L\setminus P\}=\bigcap\{M\in
MinSpe(L)|\,\, M\subseteq P\}$.\vskip 2mm
\end{center}

\noindent{\bf Proof. }The second equation is clear. It suffices to
show the first equation.

Clearly, $\bigcup\{a^{\perp}|\,\, a\in X\}\subseteq\bigcap\{P\in
Spe(L)|\,\, P\cap X=\emptyset\}$. If $x\not\in
\bigcup\{a^{\perp}|\,\, a\in X\}$ then $x\not\in a^{\perp}$ for
any $a\in X$, i.e., $x\wedge a>0$ for any $a\in X$. Consider the
set\vskip 2mm
\begin{center}
$\overline{X}=\{x\wedge a|\,\, a\in X\cup\{x\}\}$. \vskip 2mm
\end{center}
\noindent A direct computation shows that $\overline{X}$ is a
$\wedge$-semilattice with $0\not\in \overline{X}$. Let $F$ be the
filter of $L$ generated by $\overline{X}$. Then there exists an
ultrafilter $U$ of $L$ such that $U\supseteq F$. By Lemma 4.1,
$P=L\setminus U\in MinSpe(L)$ and $P\cap X=\emptyset$. Since $x\in
U$, we get $x\not\in P$, so that $x\not\in \bigcap\{P\in
Spe(L)|\,\, P\cap X=\emptyset\}$. So $\bigcup\{a^{\perp}|\,\, a\in
X\}=\bigcap\{P\in Spe(L)|\,\, P\cap X=\emptyset\}$.

Using the above results, the remains are clear. $\hfill\Box$\vskip
2mm

Now we can apply Lemma 4.1 and Lemma 4.2 to establish
characterizations of minimal prime ideals of a decomposable
lattice.\vskip 2mm

\noindent{\bf Theorem 4.3. }Let $L\in \mathbb{DL}$ and $P\in
Spe(L)$. The following conditions are equivalent:

(1) $P\in MinSpe(L)$.

(2) $P=\bigcup\{a^{\perp}|\,\, a\not\in P\}$.

(3) For any $x\in P$, $x^{\perp}\not\subseteq P$.\vskip 2mm

\noindent{\bf Proof. }(1)$\Rightarrow$(2) By Lemma 4.2, we
have\vskip 2mm
\begin{center}
$\bigcup\{a^{\perp}|\,\, a\not\in P\}=\bigcap\{M\in MinSpe(L)|\,\,
M\subseteq P\}$.\vskip 2mm
\end{center}
\noindent Since $P\in MinSpe(L)$, this means that the set $\{M\in
MinSpe(L)|\,\, M\subseteq P\}=\{P\}$, so
$P=\bigcup\{a^{\perp}|\,\, a\not\in P\}$.

(2)$\Rightarrow$(3) By (2),\vskip 2mm
\begin{center}
$P=\bigcup\{a^{\perp}|\,\, a\not\in P\}$.\vskip 2mm
\end{center}
\noindent So, for any $x\in P$, there exists some $a\not\in P$
such that $x\in a^{\perp}$. Then $a\in x^{\perp}$, which implies
$a\in x^{\perp}\setminus P$. Therefore $x^{\perp}\not\subseteq P$.

(3)$\Rightarrow$(1) Assume that $P\not\in MinSpe(L)$. Then there
exists some $M\in MinSpe(L)$ such that $P\supset M$. Pick $x\in
P\setminus M$. Then for any $0<y\in x^{\perp}$, $x\wedge y=0\in
M$. Since $M$ is prime and $x\not\in M$, we get $y\in M$, and
hence $x^{\perp}\subseteq M\subset P$, a contradiction. So $P\in
MinSpe(L)$. $\hfill\Box$\vskip 2mm

We now apply Theorem 4.3 to investigate the relationship among
prime ideals, minimal prime ideals and regular ideals. In order to
do this, we need the following two lemmas.\vskip 2mm

\noindent{\bf Lemma 4.4. }Let $L\in \mathbb{DL}$ and $0\neq A\in
Ide(L)$. Then\vskip 2mm
\begin{center}
$A^{\perp}=\bigcap\{M\in Spe(L)|\,\, A\not\subseteq
M\}=\bigcap\{P\in MinSpe(L)|\,\, A\not\subseteq P\}$.\vskip 2mm
\end{center}

\noindent{\bf Proof. }It suffices to show the first equation.

If $A\not\subseteq M$, then pick $a\in A\setminus M$, so that
$a^{\perp}\subseteq M$ since $M\in Spe(L)$. Hence
$A^{\perp}\subseteq a^{\perp}\subseteq M$. So
$A^{\perp}\subseteq\bigcap\{M\in Spe(L)|\,\, A\not\subseteq M\}$.

Now, suppose that $A^{\perp}\subset\bigcap\{M\in Spe(L)|\,\,
A\not\subseteq M\}$. Pick $0<b\in (\bigcap\{M\in Spe(L)|\,\,
A\not\subseteq M\})\setminus A^{\perp}$. Then there exists some
$0<c\in A$ such that $b\wedge c>0$. Thus $b\wedge c\in A$. Now,
Let $M\in Val(b\wedge c)$. Then $b\wedge c\not\in M$. So $b\wedge
c\not\in \bigcap\{M\in Spe(L)|\,\, A\not\subseteq M\}$, which
contradicts the fact that $b\in \bigcap\{M\in Spe(L)|\,\,
A\not\subseteq M\}$. So $A^{\perp}=\bigcap\{M\in Spe(L)|\,\,
A\not\subseteq M\}$. $\hfill\Box$\vskip 2mm

\noindent{\bf Lemma 4.5. }Let $L\in \mathbb{DL}$ and $a,b\in
L\setminus \{0\}$. The following conditions are equivalent:

(1) $a$ and $b$ are disjoint, i.e., $a\wedge b=0$.

(2) $Val(a)\cap Val(b)=\emptyset$ and $Val(a)\cup Val(b)=Val(a\vee
b)$.\vskip 2mm

\noindent{\bf Proof. }(1)$\Rightarrow$(2) Suppose that $Val(a)\cap
Val(b)\neq\emptyset$. Then there exists $M\in Val(a)\cap Val(b)$
such that $a\not\in M$ and $b\not\in M$. So $a\wedge b\not\in M$
implies $a\wedge b\neq 0$, a contradiction. Now, if $Q\in
Val(a\vee b)$ then $a\vee b\not\in Q$. Since $Q$ is an ideal of
$L$, we get that either $a\not\in Q$ or $b\not\in Q$. Without loss
of generality, assume that $a\not\in Q$. Then there exists some
$Q_{a}\in Val(a)$ such that $Q\subseteq Q_{a}$. If $Q\subset
Q_{a}$, then $a\vee b\in Q^{*}\subseteq Q_{a}$ ($Q^{*}$ denotes
the cover of $Q$ in $Ide(L)$), so that $a\in Q_{a}$, a
contradiction. So $Q=Q_{a}\in Val(a)$. Conversely, if $K\in
Val(a)\cup Val(b)$ then either $K\in Val(a)$ or $K\in Val(b)$.
Without loss of generality, assume that $K\in Val(a)$, then
$a\not\in K$, and hence $a\vee b\not\in K$. So there exists some
$Q\in Val(a\vee b)$ such that $K\subseteq Q$. If $K\subset Q$,
then $a\in K^{*}\subseteq Q$. Since $b\in P\subset Q$, $a\vee b\in
Q$, a contradiction. Therefore $P=Q\in Val(a\vee b)$.

(2)$\Rightarrow$(1) Suppose that $a\wedge b\neq 0$. Let $M\in
Val(a\wedge b)$. Then $a\not\in M$. So there exists some $P\in
Val(a)$ such that $P\supseteq M$. Similarly, $b\not\in M$. So
there exists some $Q\in Val(a)$ such that $Q\supseteq M$. By
Corollary 3.3, $P$ and $Q$ are comparable. Again, $Val(a)\cup
Val(b)=Val(a\vee b)$, so that $P=Q$, which contradicts $Val(a)\cap
Val(b)=\emptyset$. Therefore $a$ and $b$ are disjoint.
$\hfill\Box$\vskip 2mm

By induction on $n$, one can obtain that if
$\{a_{1},a_{2},\cdots,a_{n}\}$ is a mutually disjoint subset of
$L$ then
$Val(\bigvee\limits_{i=1}^{n}a_{i})=\bigcup\limits_{i=1}^{n}Val(a_{i})$.\vskip
2mm

\noindent{\bf Theorem 4.6. }Let $L\in \mathbb{DL}$ and $0\neq I\in
Ide(L)$. The following conditions are equivalent:

(1) $I$ is totally ordered.

(2) For any $0<a\in I$, $a^{\perp}=I^{\perp}$.

(3) $I^{\perp}\in Spe(L)$.

(4) $I^{\perp}\in MinSpe(L)$.

(5) $I^{\perp\perp}$ is a maximal totally ordered ideal of $L$.

(6) $I^{\perp\perp}$ is a minimal polar ideal of $L$.

(7) $I^{\perp}$ is a maximal polar ideal of $L$.

(8) For any $0<a\in I$, $a$ is special.\vskip 2mm

\noindent{\bf Proof. }(1)$\Rightarrow$(2) For any $0<a\in I$,
$a^{\perp}\supseteq I^{\perp}$ is clear. Assume that
$a^{\perp}\supset I^{\perp}$. Pick $0<x\in a^{\perp}\setminus
I^{\perp}$. Then $x\wedge a=0$ and $x\wedge b>0$ for some $b\in
I$. So $(x\wedge b)\wedge a=(x\wedge a)\wedge b=0$. On the other
hand, $0<a,x\wedge b\in I$, and hence $a$ and $x\wedge b$ are
comparable, so that\vskip 2mm
\begin{center}
$(x\wedge b)\wedge a=\min\{x\wedge b, a\}>0$.\vskip 2mm
\end{center}
This is impossible. So $a^{\perp}=I^{\perp}$.

(2)$\Rightarrow$(3) By Theorem 3.1, it suffices to show that if
$a,b\not\in I^{\perp}$ then $a\wedge b\not\in I^{\perp}$. Since
$a\not\in I^{\perp}$, there exists $0<x\in I$ such that $a\wedge
x>0$. Similarly, $b\not\in I^{\perp}$, there exists $0<y\in I$
such that and $b\wedge y>0$. We claim that $(a\wedge
x)\wedge(b\wedge y)>0$. Otherwise, $(a\wedge x)\wedge(b\wedge
y)=0$, hence $b\wedge y\in (a\wedge x)^{\perp}=I^{\perp}$ by (2),
so that $b\wedge y\in I\cap I^{\perp}=0$, a contradiction.
Therefore $I^{\perp}\in Spe(L)$.

(3)$\Rightarrow$(4) By Lemma 4.4, we have\vskip 2mm
\begin{center}
$I^{\perp}=\bigcap\{P\in MinSpe(L)|\,\, I\not\subseteq P\}$.\vskip
2mm
\end{center}
\noindent Assume that $I^{\perp}\not\in MinSpe(L)$. Then there
exists some $P\in MinSpe(L)$ such that $I^{\perp}\supset P$, so
that $I\subseteq P\subseteq I^{\perp}$. Thus $I=0$, a
contradiction. Thus $I^{\perp}\in MinSpe(L)$.

(4)$\Rightarrow$(5) We first show that $I^{\perp\perp}$ is totally
ordered. Assume that there exist $0<a,b\in I^{\perp\perp}$ such
that $a\wedge b=0$. Since $I^{\perp}$ is prime, either $a\in
I^{\perp}$ or $b\in I^{\perp}$, so that either $a=0$ or $b=0$, a
contradiction.

We next show that $I^{\perp\perp}$ is maximal. Let $J$ be a
totally ordered ideal of $L$ such that $J\supset I^{\perp\perp}$.
Pick $0<x\in J\setminus I^{\perp\perp}$. Then there exists some
$0<y\in I^{\perp}$ such that $x\wedge y>0$. Now, pick $0<a\in I$.
Then $(x\wedge y)\wedge a=0$ since $x\wedge y\in I^{\perp}$. On
the other hand, $0<x\wedge y\in J, a\in I\subseteq
I^{\perp\perp}\subseteq J$ and $J$ is totally ordered, so that
$(x\wedge y)\wedge a=\min\{x\wedge y,a\}>0$. This is impossible.
Therefore $I^{\perp\perp}$ is a maximal totally ordered ideal of
$L$.

(5)$\Rightarrow$(6) Let $D\in P(L)$ be such that $D\subseteq
I^{\perp\perp}$. Then $D$ is totally ordered. By using the result
of
(1)$\Rightarrow$(2)$\Rightarrow$(3)$\Rightarrow$(4)$\Rightarrow$(5),
$D=D^{\perp\perp}$ is a maximal totally ordered ideal of $L$, so
that $D=I^{\perp\perp}$. So $I^{\perp\perp}$ is a minimal polar
ideal of $L$.

(6)$\Rightarrow$(7) Since the map $P\rightarrow P^{\perp}$ for any
$P\in P(L)$ is a dual isomorphism of lattices, $I^{\perp\perp}$ is
a minimal polar ideal of $L$ implies that $I^{\perp}$ is a maximal
polar ideal of $L$.

(7)$\Rightarrow$(8) For any $0<a\in I$, assume that $a$ has two
distinct values $Q_{1}$ and $Q_{2}$. Since $a\not\in Q_{1}$ and
$a\wedge b=0$ for any $b\in I^{\perp}$, so that
$I^{\perp}\subseteq Q_{1}$. Similarly, $I^{\perp}\subseteq Q_{2}$.
Since $Q_{1}$ and $Q_{2}$ are incomparable, we may pick\vskip 2mm
\begin{center}
$0<x\in Q_{1}\setminus Q_{2}, 0<y\in Q_{2}\setminus Q_{1}$ with
$x\wedge y=0$.\vskip 2mm
\end{center}
\noindent So $x^{\perp}=y^{\perp}=I^{\perp}$ by the maximality of
$I^{\perp}$. Again, $x\wedge y=0$ implies that $x,y\in
I^{\perp}\subseteq Q_{1}\cap Q_{2}$, a contradiction. So $a$ is
special.

(8)$\Rightarrow$(1) Assume that $I$ is not totally ordered. Then
there exist $0<a,b\in I$ such that $a\wedge b=0$. By Lemma 4.5,
$Val(a\vee b)=Val(a)\cup Val(b)$, i.e., $a\vee b$ has at least two
distinct values, a contradiction. Therefore $I$ must be totally
ordered. $\hfill\Box$\vskip 2mm

By using Theorem 4.6, we shall investigate the relationship
between polar ideals and minimal prime ideals of a decomposable
lattice.\vskip 2mm

\noindent{\bf Theorem 4.7. }Let $L\in \mathbb{DL}$. If for any
$P,Q\in P(L)$ either $L=P\vee Q$ or $P$ and $Q$ are comparable,
then every polar ideal of $L$ is minimal prime, i.e.,
$P(L)\subseteq MinSpe(L)$.\vskip 2mm

\noindent{\bf Proof. }By way of contradiction. If there exists
$P\in P(L)$ such that $P\not\in MinSpe(L)$, write $P=A^{\perp}$,
then $A$ is not totally ordered by Theorem 4.6. So there exist
$0<a,b\in A$ such that $a\wedge b=0$. We divide the proof into two
steps.

{\bf Step 1. }If $a^{\perp}$ and $b^{\perp}$ are incomparable then
$L=a^{\perp}\vee b^{\perp}$ by hypothesis. Clearly, $a^{\perp}$
and $a^{\perp\perp}$ are incomparable, then $L=a^{\perp}\vee
a^{\perp\perp}$. So\vskip 2mm
\begin{center}
$a^{\perp\perp}=a^{\perp\perp}\cap
L=a^{\perp\perp}\cap(a^{\perp}\vee b^{\perp})=a^{\perp\perp}\cap
b^{\perp}\subseteq b^{\perp}$.\vskip 2mm
\end{center}
\noindent Similarly, $b^{\perp\perp}\subseteq a^{\perp}$.\vskip
2mm

(i) If $a^{\perp\perp}$ and $b^{\perp\perp}$ are incomparable then
$L=a^{\perp\perp}\vee b^{\perp\perp}$. So\vskip 2mm
\begin{center}
$a^{\perp}=a^{\perp}\cap L=a^{\perp\perp}\cap(a^{\perp\perp}\vee
b^{\perp\perp})=a^{\perp}\cap b^{\perp\perp}\subseteq
b^{\perp\perp}$.\vskip 2mm
\end{center}
\noindent Thus $a^{\perp}=b^{\perp\perp}$. Similarly,
$a^{\perp\perp}=b^{\perp}$. It follows that $a^{\perp}\cap
b^{\perp}=a^{\perp}\cap a^{\perp\perp}=\{0\}$. So
$A^{\perp}\subseteq a^{\perp}\cap b^{\perp}=\{0\}$, a
contradiction.\vskip 2mm

(ii) If $a^{\perp\perp}$ and $b^{\perp\perp}$ are comparable then
$a^{\perp\perp}\subseteq b^{\perp\perp}=a^{\perp}$, and hence
$a^{\perp\perp}=0$ or $b^{\perp\perp}\subseteq
a^{\perp\perp}=b^{\perp}$, and hence $b^{\perp\perp}=0$. It
follows that $a^{\perp}=L$ or $b^{\perp}=L$, this is
impossible.\vskip 2mm

{\bf Step 2. }If $a^{\perp}$ and $b^{\perp}$ are comparable then
$a^{\perp}\subseteq b^{\perp}$ or $b^{\perp}\subseteq a^{\perp}$.
So $b\in a^{\perp}\subseteq b^{\perp}$ and hence $b=0$ or $a\in
b^{\perp}\subseteq a^{\perp}$ and hence $a=0$, this is also
impossible.

In view of Step 1 and Step 2, $A$ is totally ordered. So $P\in
MinSpe(L)$.$\hfill\Box$ \vskip 2mm

Recall that a lattice $L$ is called projectable if
$L=x^{\perp}\vee x^{\perp\perp}$ for any $x\in L$. We denote by
$\mathbb{T}$ the class of projectable lattices. \vskip 2mm

\noindent{\bf Theorem 4.8. }Let $L\in \mathbb{DL}$. The following
conditions are equivalent:

(1) $Spe(L)=MinSpe(L)$.

(2) $L=(a]\vee a^{\perp}$ for any $a\in L$.

(3) $L\in \mathbb{T}$ and $(x]=x^{\perp\perp}$ for any $x\in L$.
\vskip 2mm

\noindent{\bf Proof. }(1)$\Rightarrow$(2) Assume that there exists
$a\in L$ such that $(a]\vee a^{\perp}\subset L$. Pick $x\in
L\setminus(a]\vee a^{\perp}$. Then there exists $M\in Val(x)$ such
that $(a]\vee a^{\perp}\subseteq M$. By (1), $M$ is minimal prime.
But $a\in M$ and $a^{\perp}\subseteq M$, which contradicts Theorem
4.3.

(2)$\Rightarrow$(3) Clearly, $L\in \mathbb{T}$. Now, given any
$x\in L$, $(x]\subseteq x^{\perp\perp}$ is clear. Again,\vskip 2mm
\begin{center}
$x^{\perp\perp}=x^{\perp\perp}\cap L=x^{\perp\perp}\cap((x]\vee
x^{\perp})=(x^{\perp\perp}\cap(x])\vee (x^{\perp\perp}\cap
x^{\perp})=x^{\perp\perp}\cap(x]\subseteq(x]$,\vskip 2mm
\end{center}
\noindent so $x^{\perp\perp}=(x]$.

(3)$\Rightarrow$(1) By (3), $L=(x]\vee x^{\perp}$ for any $x\in
L$. Assume that there exists $P\in Spe(L)$ such that $P$ is not
minimal. Then there exists some $M\in MinSpe(L)$ such that
$P\supset M$. Pick $a\in P\setminus M$. Then $a^{\perp}\subseteq
M\subset P$, so that $L=(a]\vee a^{\perp}\subseteq P$, a
contradiction. Therefore $Spe(L)=MinSpe(L)$.$\hfill\Box$ \vskip
2mm

Recall that a minimal element of a partially ordered set is an
atom. If every element exceeds an atom, the partially ordered set
is called atomic. Theorem 3.1 shows that the set $Spe(L)$ of all
prime ideals of a decomposable lattice $L$ is an atomic root
system under inclusion. It is natural to ask under what condition
to make $V(L)$ atomic, i.e., every regular ideal of $L$ contains a
minimal regular ideal. In order to do this, we need the following
lemma. Since its proof is direct, we shall omit it.\vskip 2mm

\noindent{\bf Lemma 4.9. }Let $L\in \mathbb{DL}$ and $M\in
Ide(L)$. Then $M\in MinSpe(L)$ if and only if there exists a
maximal chain $\{M_{\lambda}\}_{\lambda\in\Lambda}$ of $V(L)$ such
that $M=\bigcap\limits_{\lambda\in\Lambda}M_{\lambda}$. \vskip 2mm

\noindent{\bf Theorem 4.10. }Let $L\in \mathbb{DL}$. If each prime
ideal in $L$ contains a finite number of minimal prime ideals then
the following conditions are equivalent:

(1) Every minimal prime ideal of $L$ is regular, i.e.,
$MinSpe(L)\subseteq V(L)$.

(2) $V(L)$ is atomic. \vskip 2mm

\noindent{\bf Proof. }(1)$\Rightarrow$(2) is clear.

(2)$\Rightarrow$(1) Given any $P\in MinSpe(L)$, by Lemma 4.9,
there exists a maximal chain of $V(L)$, write $\{Q_{\gamma}\in
V(L)|\,\, \gamma\in \bigtriangleup\}$, such that
$P=\bigcap\limits_{\gamma\in\bigtriangleup}Q_{\gamma}$. Pick
$Q_{\gamma_{1}}\in \bigtriangleup$. Since $V(L)$ is atomic,
$Q_{\gamma_{1}}$ contains an atom, write $Q_{1}$. If
$Q_{\gamma}\supseteq Q_{1}$ for any $\gamma\in\bigtriangleup$,
then $P=Q_{1}$, we are done. Otherwise, there exists
$Q_{\gamma_{2}}\in \bigtriangleup$ such that
$Q_{\gamma_{2}}\subset Q_{\gamma_{1}}$, but $Q_{1}\not\subseteq
Q_{\gamma_{2}}$. Similarly, $Q_{\gamma_{2}}$ contains an atom,
write $Q_{2}$. If $Q_{\gamma}\supseteq Q_{2}$ for any
$\gamma\in\bigtriangleup$, then $P=Q_{2}$, we are done. We claim
that this process must end. Otherwise, we may obtain an infinite
number of atoms in $V(L)$, write
$\{Q_{1},Q_{2},\cdots,Q_{n},\cdots\}$ which satisfy $Q_{i}\neq
Q_{j}$ for any $i\neq j$. Clearly, $Q_{\gamma_{1}}$ contains each
$Q_{i}$ for all $i=1,2,\cdots,n,\cdots$, a contradiction. So each
minimal prime ideal in $L$ is regular. $\hfill\Box$\vskip 8mm

\noindent{\bf 5. Special ideals}\vskip 4mm

In this section, we characterize special ideals of a decomposable
lattice and then investigate the relationship among prime ideals,
minimal prime ideals, regular ideals and special ideals.

Recall that for a lattice $L$ and $0<x\in L$, if $M$ is the unique
value of $x$, $M$ or $x$ is called special. We denote by $S(L)$
the set of all special ideals of the lattice $L$.\vskip 2mm

\noindent{\bf Theorem 5.1. }Let $L$ be a lattice and $M\in
Ide(L)$. Then the following conditions are equivalent:

(1) $M\in S(L)$.

(2) If $\bigcap\limits_{\lambda\in\Lambda}I_{\lambda}\subseteq M$,
where $\{I_{\lambda}\}_{\lambda\in\Lambda}\subseteq Ide(L)$, then
$I_{\lambda}\subseteq M$ for some $\lambda$.

(3) $M$ is the unique value of $x$, where $x\in M^{*}\setminus
M$.\vskip 2mm

\noindent{\bf Proof. }(1)$\Leftrightarrow$(2) is clear. It
suffices to show (1)$\Leftrightarrow$(3)

(1)$\Rightarrow$(3) Consider the set\vskip 2mm
\begin{center}
$\bigtriangleup=\{J_{\lambda}\in Ide(L)|\,\
J_{\lambda}\not\subseteq M, \,\, \lambda\in \Lambda\}$.\vskip 2mm
\end{center}
\noindent Since $M$ is special, $0\neq {\bigcap\limits_{\lambda\in
\Lambda}} J_{\lambda}\not\subseteq M$, and hence pick $0\neq x\in
({\bigcap\limits_{\lambda\in \Lambda}} J_{\lambda})\setminus M$.
Now, if $K$ is an ideal of $L$ with respect to not containing $x$,
and $K\supset M$, then $K\in \bigtriangleup$, so that $x\in K$, a
contradiction. So $M$ is a maximal ideal with respect to not
containing $x$. Again, if $N$ is any maximal ideal of $L$ with
respect to not containing $x$, then since $x\not\in N$, $N\not\in
\bigtriangleup$. So $N\subseteq M$, and hence $N=M$. Therefore $M$
is the unique maximal ideal of $L$ with respect to not containing
$x$ and clearly $x\in M^{*}\setminus M$.

(3)$\Rightarrow$(1) Clearly, $M$ is regular. Now, let
$\{I_{\lambda}\}_{\lambda\in \Lambda}$ be any nonempty family of
ideals of $L$ such that ${\bigcap\limits_{\lambda\in \Lambda}}
I_{\lambda} \subseteq M$. Since $x\not\in M$, there exists some
$\lambda\in \Lambda$ such that $x\not\in I_{\lambda}$. So there
exists some $N\in Val(x)$ such that $I_{\lambda} \subseteq N$. By
assumption, $M$ is the unique maximal ideal with respect to not
containing $x$, so that $N=M$. Therefore $I_{\lambda} \subseteq
M$. So $M$ is special.$\hfill\Box$\vskip 2mm

In order to investigate the relationship among prime ideals,
minimal prime ideals and special ideals of a decomposable lattice,
we need the following two lemmas.

For a lattice $L$ and $P\in Spe(L)$, write $S_{P}=\bigcap\{M\in
MinSpe(L)|\,\, M\subseteq P\}$.\vskip 2mm

\noindent{\bf Lemma 5.2. }Let $L\in \mathbb{DL}$ and
$P_{1},P_{2}\in Spe(L)$. Then $S_{P_{1}}\subseteq P_{2}$ if and
only if $P_{1}$ and $P_{2}$ are comparable.\vskip 2mm

\noindent{\bf Proof. }The sufficiency is clear. For the necessity,
assume that $P_{1}\parallel P_{2}$. Pick $a_{1}\in P_{2}\setminus
P_{1}, a_{2}\in P_{1}\setminus P_{2}$ with $a_{1}\wedge a_{2}=0$.
Now, let $M\in MinSpe(L)$ be such that $M\subseteq P_{2}$. Then
$a_{1}\in M$, so $a_{1}\in S_{P_{2}}$. Similarly, $a_{2}\in
S_{P_{1}}$. Since $S_{P_{1}}\subseteq P_{2}$, we get $a_{2}\in
P_{2}$, a contradiction.$\hfill\Box$\vskip 2mm

\noindent{\bf Lemma 5.3. }Let $L\in \mathbb{DL}$ and $P\in
Spe(L)$. Then $S_{P}=\{a\in L|\,\, a=0$ or for any $Q\in Spe(L)$
with $a\not\in Q, Q$ and $P$ are not comparable $\}$.\vskip 2mm

\noindent{\bf Proof. }Write $K=\{a\in L|\,\, a=0$ or for any $Q\in
Spe(L)$ with $a\not\in Q, Q$ and $P$ are not comparable $\}$. If
$S_{P}\not\subseteq K$, pick $0<a\in S_{P}\setminus K$, then there
exists $Q\in Val(a)$ such that $Q$ and $P$ are comparable. If
$Q\subseteq P$ then $a\in S_{P}\subseteq Q$ by Lemma 5.2, a
contradiction. If $P\subset Q$ then $a\in S_{P}\subseteq P\subset
Q$, a contradiction. So $S_{P}\subseteq K$. Conversely, if
$K\not\subseteq S_{P}$, pick $0<b\in K\setminus S_{P}$, then there
exists $M\in MinSpe(L)$ with $M\subseteq P$ such that $b\not\in
M$. But $b\in K$, so that $M$ and $P$ are comparable, a
contradiction. So $S_{P}=K$.$\hfill\Box$\vskip 2mm

\noindent{\bf Theorem 5.4. }Let $L\in \mathbb{DL}$ and $I\in
Ide(L)$. The following conditions are equivalent:

(1) There exists a unique value $Q$ of $g$ such that $Q\supseteq
I$, and for any $x\in L\setminus I$, $x\wedge g\not\in I$.

(2) $S_{P}\subseteq I\subseteq P$, where $P\in Val(g)$.\vskip 2mm

\noindent{\bf Proof. }(1)$\Rightarrow$(2) Let $P$ be the unique
value of $g$ containing $I$. Since $I=\bigcap\{N\in V(L)|\,\,
I\subseteq N\}$, it suffices to show that if $N\in V(L)$ with
$N\supseteq I$ then $N\supseteq S_{P}$.

Suppose that there exists $N\in V(L)$ with $I\subseteq N$, but
$S_{P}\not\subseteq N$. Then, by Lemma 5.2, $P\parallel N$. Pick
$0<x\in N^{*}\setminus N$ and $0<y\in P\setminus N$. Then
$0<x\wedge y\in (N^{*}\setminus N)\cap P$. Using this method, we
see that there exist\vskip 2mm
\begin{center}
$0<a\in (N^{*}\setminus N)\cap P$ and $0<b\in (P^{*}\setminus
P)\cap N$.\vskip 2mm
\end{center}
\noindent Since $L\in \mathbb{DL}$, we may further assume that
$a\wedge b=0$. Now, if $a$ has a value $K$ such that $K\subseteq
P$ then since $a\not\in K$ implies $b\in K\subseteq P$, a
contradiction. So each value of $a$ is not comparable with $P$. By
Lemma 5.3, $a\in S_{P}$. So $a\in(N^{*}\setminus N)\cap S_{P}$.
Now, let $0<x\in (N^{*}\setminus N)\cap S_{P}$. By (1), $x\wedge
g\not\in I$. Then there exists $K_{x\wedge g}\in Val(x\wedge g)$
such that $K_{x\wedge g}\supseteq I$. So $g\not\in K_{x\wedge g}$,
which implies that there exists some $K_{g}\in Val(g)$ such that
$K_{g}\supseteq K_{x\wedge g}\supseteq I$. By (1), $K_{g}=P$. On
the other hand, $x\not\in K_{x\wedge g}$, there exists $K_{x}\in
Val(x)$ such that $K_{x}\supseteq K_{x\wedge g}$, which implies
that $K_{x}$ and $P$ are comparable. So $x\in S_{P}\subseteq
K_{x}$, a contradiction. So $S_{P}\subseteq I\subseteq P$.

(2)$\Rightarrow$(1)Assume that that there exists another value
$P_{1}$ of $g$ such that $S_{P}\subseteq I\subseteq P_{1}$. Note
that $P_{1}\neq P$ implies that $P\parallel P_{1}$. But, by Lemma
5.2, $S_{P}\subseteq P_{1}$ implies that $P$ and $P_{1}$ are
comparable, a contradiction. So $P=P_{1}$, and hence $P$ is the
only value of $g$ containing $I$.

Now, let $\Gamma=\{M\in V(L)|\,\, I\subseteq M\}$. Clearly
$I=\bigcap\{M\in V(L)|\,\, I\subseteq M\}$, which implies that $M$
and $P$ are comparable for any $M\in \Gamma$. Set
$\Gamma_{0}=\Gamma\setminus \{M\in V(L)|\,\, P\subset M\}$.
Clearly $I=\bigcap\limits_{M\in T_{0}}M$. Then for any $M\in
\Gamma_{0}$, $M\subseteq P$. So if $x\not\in I$ then there exists
$M\in \Gamma_{0}$ such that $x\not\in M$. Again, $g\not\in P$,
then $g\not\in M$. So $x\wedge g\not\in M$, and hence $x\wedge
g\not\in I$, as desired.$\hfill\Box$ \vskip 2mm

\noindent{\bf Theorem 5.5. }Let $L\in \mathbb{DL}$ and $K\in
Ide(L)$. The following conditions are equivalent:

(1) $K\in Spe(L)$ and for any $x\in L\setminus K$, $x>K$.

(2) $K\in Spe(L)$ and for any $I\in Ide(L)$, $K$ and $I$ are
comparable.

(3) For any $L\neq P\in P(L)$, $P\subseteq K$.

(4) For any $M\in MinSpe(L)$, $M\subseteq K$.

(5) For any $a\in L\setminus K$, $a^{\perp}=\{0\}$.

(6) For any $a\in L\setminus K$, $a$ is special.\vskip 2mm

\noindent{\bf Proof. }(1)$\Rightarrow$(2) Assume that there exists
some $I\in Ide(L)$ such that $K$ and $I$ are incomparable. Pick
$x\in I\setminus K$. By (1), $x>K$, so that $K\subseteq I$, a
contradiction.

(2)$\Rightarrow$(3) Given any $L\neq P\in P(L)$, if
$P\not\subseteq K$ then $K\subset P$ by (2). Pick $x\in P\setminus
K$. Since $K\in Spe(L)$ and $x\not\in K$, $P^{\perp}\subseteq
K\subset P$. So $P^{\perp}=0$, and hence $P=P^{\perp\perp}=L$, a
contradiction.

(3)$\Rightarrow$(4) For any $M\in MinSpe(L)$, by Theorem
4.3,\vskip 2mm
\begin{center}
$M=\bigcup\{a^{\perp}|\,\, a\not\in M\}$.\vskip 2mm
\end{center}
\noindent By (3), $M\subseteq K$.

(4)$\Rightarrow$(5) Given any $a\in L\setminus K$, by (4),
$M\subseteq K$ for any $M\in MinSpe(L)$. Then $a\not\in M$ for any
$M\in MinSpe(L)$, so that $a^{\perp}\subseteq M$ for any $M\in
MinSpe(L)$. So $a^{\perp}\subseteq\bigcap MinSpe(L)=0$, i.e.,
$a^{\perp}=\{0\}$.

(5)$\Rightarrow$(6) Given any $a\in L\setminus K$, assume that $a$
is not special. Then $a$ has at least two distinct values
$Q_{1},Q_{2}$. Clearly, $Q_{1}\parallel Q_{2}$. Pick\vskip 2mm
\begin{center}
$a_{1}\in Q_{1}\setminus Q_{2}$ and $a_{2}\in Q_{2}\setminus
Q_{1}$ with $a_{1}\wedge a_{2}=0$.\vskip 2mm
\end{center}
\noindent Clearly, $a_{1},a_{2}\not\in K$, but
$a_{1}^{\perp}\neq\{0\}$, a contradiction.

(6)$\Rightarrow$(1) We first show that $K\in Spe(L)$. Assume that
there exist $0<a,b\in L$ such that $a\wedge b=0$, but $a\not\in K$
and $b\not\in K$. Then $a\vee b\not\in K$. Notice that $a\wedge
b=0$, so that $Val(a\vee b)=Val(a)\cup Val(b)$. So $a\vee b$ is
not special, a contradiction.

We next show that for any $x\in L\setminus K$, $x>K$. Otherwise,
there exists $0<k\in K$ such that $x\parallel k$. Since $L\in
\mathbb{DL}$, we may further assume that $x\wedge k=0$. Clearly,
$x\vee k\not\in K$. But $Val(x\vee k)=Val(x)\cup Val(k)$, it
follows that $x\vee k$ is not also special, which ends the
proof.$\hfill\Box$\vskip 2mm

Recall that if $L\in \mathbb{DL}$ then the set $Ide(L)$ of all
ideals of $L$ is a distributive lattice by the rule: $I\wedge
J=I\cap J$ and $I\vee J=\{a\vee b|\,\, a\in I, b\in J\}$. So
$Ide(L)$ is $\alpha$-distributive, i.e., for any $I\in Ide(L)$ and
any subset $\{J_{\lambda}\}_{\lambda\in \Lambda}\subseteq Ide(L)$
with $|\Lambda|=\alpha$, $I\bigcap({\bigvee\limits_{\lambda\in
\Lambda }}J_{\lambda})={\bigvee\limits_{\lambda\in
\Lambda}}(I\bigcap J_{\lambda})$. But, in general, it is not dual
$\alpha$-distributive, i.e., $I\bigvee({\bigcap\limits_{\lambda\in
\Lambda }}J_{\lambda})={\bigcap\limits_{\lambda\in
\Lambda}}(I\bigvee J_{\lambda})$  does not hold.

In order to establish the condition that $V(L)=S(L)$, let us
recall that a lattice $L$ is called completely distributive if for
any nonempty family $\{a_{i,j}\}_{i\in I,j\in J}\subseteq L$,
whenever ${\bigvee\limits_{i\in I}}\,{\bigwedge\limits_{j\in
J}}a_{i,j}$ and ${\bigwedge\limits_{f\in
I^{J}}}\,{\bigvee\limits_{i\in I}}a_{i,f(i)}$ exist in $L$,
then\vskip 2mm
\begin{center}
${\bigvee\limits_{i\in I}}\,{\bigwedge\limits_{j\in
J}}a_{i,j}={\bigwedge\limits_{f\in J^{I}}}\,{\bigvee\limits_{i\in
I}}a_{i,f(i)}$,\vskip 2mm
\end{center}
\noindent where $J^{I}$ denotes the set of all maps from $I$ to
$J$.\vskip 2mm

\noindent{\bf Theorem 5.6. }Let $L\in \mathbb{DL}$. The following
conditions are equivalent:

(1) $V(L)=S(L)$.

(2) $Ide(L)$ is completely distributive.

(3) $Ide(L)$ is $\alpha$-distributive. \vskip 2mm

\noindent{\bf Proof. }(1)$\Rightarrow$(2) Let $\{K_{i,j}\}_{i\in
I,j\in J}$ be any nonempty family of ideals of $L$, and suppose
that ${\bigvee\limits_{i\in I}}\,\ {\bigcap\limits_{j\in
J}}K_{i,j}$\,\,\,\, and\,\,\,\, ${\bigcap\limits_{f\in I^{J}}}\,\
{\bigvee\limits_{i\in I}}K_{i,f(i)}$ exist in $Ide(L)$.
Write\vskip 2mm
\begin{center}
$A={\bigvee\limits_{i\in I}}\,\ {\bigcap\limits_{j\in J}}K_{i,j}$,
and $B={\bigcap\limits_{f\in I^{J}}}\,\ {\bigvee\limits_{i\in
I}}K_{i,f(i)}$.\vskip 2mm
\end{center}
\noindent Clearly, $A\subseteq B$. Since for any $I\in Ide(L)$,
$I=\bigcap\{M\in V(L)|\,\, I\subseteq M \}$ and thus it suffices
to show that for any $M\in V(L)$, if $M\supseteq A$ then
$M\supseteq B$. Now, suppose $A\subseteq M$; then
${\bigcap\limits_{j\in J}}K_{i,j}\subseteq M$ for any $i\in I$. By
assumption, $M\in V(L)=S(L)$, so there exists some $j_{i}\in J$
such that $K_{i,j_{i}}\subseteq M$. Now let $f(i)=j_{i}$ for any
$i\in I$; then $K_{i,f(i)}\subseteq M$ for any $i\in I$. It
follows that ${\bigvee\limits_{i\in I}}K_{i,f(i)}\subseteq M$. So
we get $B\subseteq M$. Consequently, we obtain $A=B$.

(2)$\Rightarrow$(3) is clear.

(3)$\Rightarrow$(1) Let $M\in V(L)$ and let
$\{I_{\lambda}\}_{\lambda\in \Lambda}$ be any nonempty family of
ideals of $L$ with $|\Lambda|=\alpha$ such that
${\bigcap\limits_{\lambda\in \Lambda}} I_{\lambda} \subseteq M$.
Since $R$ is dual $\alpha$-distributive, we then have\vskip 2mm
\begin{center}
$M=M\bigvee({\bigcap\limits_{\lambda\in \Lambda}}
I_{\lambda})={\bigcap\limits_{\lambda\in \Lambda}} (M\bigvee
I_{\lambda})$.\vskip 2mm
\end{center}
\noindent So there exists some $\lambda_{0}\in \Lambda$ such that
$M\vee I_{\lambda_{0}}=M$, i.e., $I_{\lambda_{0}}\subseteq M$. So
$M\in S(L)$. Therefore $V(L)=S(L)$.$\hfill\Box$\vskip 2mm

At the end of this paper, we shall investigate decomposable
lattices in which each nonzero element has only finitely many
values.\vskip 2mm

\noindent{\bf Lemma 5.7. }Let $L\in \mathbb{DL}$. If
$Q_{1},Q_{2},\cdots, Q_{n}$ are mutually incomparable prime ideals
of $L$ and $a\not\in Q_{i}$ for $i=1,2,\cdots,n$, then there exist
$a_{i}\in (\bigcap\limits_{j\neq i}Q_{j}) \setminus Q_{i}$ such
that $0<a_{i}<a$ for $i=1,2,\cdots,n$ and $a_{i}\wedge a_{j}=0$
for $i\neq j$.\vskip 2mm

\noindent{\bf Proof. }By induction on $n$. If $n=2$ then pick
$0<x_{1}\in Q_{2}\setminus Q_{1}, 0<x_{2}\in Q_{1}\setminus
Q_{2}$. Clearly, $x_{1}\parallel x_{2}$, so there exist
$y_{1},y_{2}\in L$ such that\vskip 2mm
\begin{center}
$x_{1}=y_{1}\vee(x_{1}\wedge x_{2}), x_{2}=y_{2}\vee(x_{1}\wedge
x_{2})$ and $y_{1}\wedge y_{2}=0$.\vskip 2mm
\end{center}
\noindent Now, set $a_{i}=a\wedge y_{i}$ for $i=1,2$. Then
$0<a_{1}\in Q_{2}\setminus Q_{1}, 0<a_{2}\in Q_{1}\setminus Q_{2}$
with $0<a_{i}<a$ and $a_{1}\wedge a_{2}=0$ for $i=1,2$.

Assume that the conclusion holds for the case $n-1$. Now consider
the case $n$. We divide the proof into three steps.\vskip 2mm

{\bf Step 1. }For prime ideals $Q_{1},Q_{2},\cdots, Q_{n-1}$,
there exist $b_{i}\in (\bigcap\limits_{1\leq j\neq i\leq
n-1}Q_{j})\setminus Q_{i}$ such that $0<b_{i}<a$ for
($i=1,2,\cdots,n-1$) and $b_{i}\wedge b_{j}=0$ for $i\neq
j$.\vskip 3mm

{\bf Step 2. }For prime ideals $Q_{2},Q_{3},\cdots, Q_{n}$, there
exist $c_{i}\in (\bigcap\limits_{2\leq j\neq i\leq
n}Q_{j})\setminus Q_{i}$ such that $0<c_{i}<a$ for
$i=2,3,\cdots,n$ and $c_{i}\wedge c_{j}=0$ for $i\neq j$. \vskip
3mm

{\bf Step 3. }Set $a_{i}=b_{i}\wedge c_{i}$ for
$i=2,3,\cdots,n-1$. Clearly, $a_{i}\in (\bigcap\limits_{1\leq
j\neq i\leq n}Q_{j})\setminus Q_{i}$ with $0<a_{i}<a$ for
$i=2,3,\cdots,n-1$ and $a_{i}\wedge a_{j}=0$ for $i\neq j$.\vskip
3mm

Last, for prime ideals $Q_{1}, Q_{n}$, since $Q_{1}\parallel
Q_{n}$, pick $0<f_{1}\in Q_{n}\setminus Q_{1}, 0<f_{n}\in
Q_{1}\setminus Q_{n}$ with $f_{1}\wedge f_{n}=0$. Set\vskip 2mm
\begin{center}
$a_{1}=f_{1}\wedge b_{1},\,\,\,\, a_{n}=f_{n}\wedge c_{n}$.\vskip
2mm
\end{center}
\noindent Then $a_{i}\in (\bigcap\limits_{j\neq i}Q_{j}) \setminus
Q_{i}$ with $0<a_{i}<a$ for $i=1,2,\cdots,n$ and $a_{i}\wedge
a_{j}=0$ for $i\neq j$, which completes the
proof.$\hfill\Box$\vskip 2mm

\noindent{\bf Lemma 5.8. }Let $L\in \mathbb{DL}$ and $0<a\in L$.
If $a$ has only $n$ values $Q_{1},Q_{2},\cdots, Q_{n}$ then
$a=\bigvee\limits_{i=1}^{n}a_{i}$ and each $Q_{i}$ is the only
value of $a_{i}$ for $i=1,2,\cdots,n$ and $a_{i}\wedge a_{j}=0$
for $i\neq j$.\vskip 2mm

\noindent{\bf Proof. }Clearly, $Q_{1},Q_{2},\cdots, Q_{n}$ are
mutually incomparable prime ideals of $L$ and $a\not\in Q_{i}$ for
$i=1,2,\cdots,n$. By Lemma 5.7, there exist $a_{i}\in
(\bigcap\limits_{j\neq i}Q_{j}) \setminus Q_{i}$ such that
$0<a_{i}<a$ for $i=1,2,\cdots,n$ and $a_{i}\wedge a_{j}=0$ for
$i\neq j$. Clearly, each $Q_{i}$ is a value of $a_{i}$ for
$i=1,2,\cdots,n$. Assume that $a_{i}$ has another value, write
$Q_{0}$. Then $Q_{0}\in Val(a)$, so there exists some $Q_{j}\in
Val(a)$ with $j\neq i$ such that $Q_{0}\subseteq Q_{j}$, and hence
$a_{j}\in Q_{0}\subseteq Q_{j}$, a contradiction.

Finally, we show that $a=\bigvee\limits_{i=1}^{n}a_{i}$. Since
$0<a_{i}<a$ for $i=1,2,\cdots,n$, we have
$\bigvee\limits_{i=1}^{n}a_{i}\leq a$. Assume that
$\bigvee\limits_{i=1}^{n}a_{i}<a$. Then $a\not\in
(\bigvee\limits_{i=1}^{n}a_{i}]$. So there exists some $Q_{i}\in
Val(a)$ such that $(\bigvee\limits_{i=1}^{n}a_{i}]\subseteq
Q_{i}$, and hence $a_{i}\in
(\bigvee\limits_{i=1}^{n}a_{i}]\subseteq Q_{i}$, a contradiction.
So $a=\bigvee\limits_{i=1}^{n}a_{i}$.$\hfill\Box$ \vskip 2mm

By Lemma 5.7 and Lemma 5.8, we have\vskip 2mm

\noindent{\bf Theorem 5.9. }Let $L\in \mathbb{DL}$. Then the
following conditions are equivalent:

(1) Each nonzero element in $L$ has only finitely many values.

(2) For any $0<a\in L$, $a=a_{1}\vee a_{2}\vee\cdots\vee a_{n}$,
where $a_{i}\wedge a_{j}=0$ for $i\neq j$ and each $a_{i}$ is
special.\vskip 8mm

\noindent{\bf Acknowledgement }\vskip 4mm

This research was supported by the Science and Technology
Development Foundation of Nanjing University of Science and
Technology and the National Natural Science Foundation of China.
\vskip 8mm

\baselineskip 11pt \vskip 0.88 true cm

\end{document}